\documentclass[a4paper, oneside,12pt]{article}
\title{Classification of Conic Sections in $PE_2(\mathbb{R})$}
\author{\and J.
Beban-Brki\'{c} \and M. \v{S}imi\'{c} Horvath}
\date{}

\usepackage{latexsym}
\usepackage{amssymb}
\usepackage{amsmath}

\usepackage[pdftex]{graphicx}

\usepackage{fancyhdr}
\newtheorem{defi}{Definition}
\newtheorem{tm}{Theorem}

\newtheorem{prop}{Proposition}

\def\naslov#1{
    \if@openright\cleardoublepage\else\clearpage\fi
    \thispagestyle{plain}
    \addcontentsline{toc}{chapter}{\protect \numberline{#1}}
    {
        \vspace*{50pt}
        \parindent 0pt
        \Huge \bfseries #1 \par\nobreak
        \vskip 40pt
        \normalfont
        \normalsize
    }
}

\begin{document}
\maketitle \textbf{Abstract:} {\small
This paper gives a complete classification of conics in $PE_2(\mathbb{R})$. The classification has been made
earlier (Reveruk \cite{rew}), but it showed to be
incomplete and not possible to cite and use in
further studies of properties of conics, pencil of conics, and of
quadratic forms in pseudo-Euclidean spaces. This paper provides that.
A pseudo-orthogonal matrix, pseudo-Euclidean values of a
matrix, diagonalization of a matrix in a pseudo-Euclidean way are
introduced. Conics are divided in families and by
types, giving both of them geometrical meaning. The invariants of a conic with respect to the group of
motions in $PE_2(\mathbb{R})$ are determined, making it possible to determine a
conic without reducing its equation to canonical form. An overview table is given.\\
 \vspace{4mm}\textbf{Key words:} {\small pseudo-Euclidean plane $PE_2(\mathbb{R})$,
conic section}\\ \vspace{2.5 mm} \textbf{MSC 2010:}{\small 51A05,
51N25 }

\newpage
\section{Pseudo-Euclidean plane}
The pseudo-Euclidean plane is a real affine plane where the metric
is introduced by the \textit{absolute figure} $(\omega, \Omega _1, \Omega_2)$ consisting of the
line $\omega$ at infinity and the points $\Omega_1, \Omega_2 \in
\omega$. Any line passing through $\Omega _1$ or $\Omega_2$ is called an \textit{isotropic} line and any point incident with $\omega$ is called an \textit{isotropic} point.\\
Let $T=(x_0:x_1:x_2)$ denote any point in the plane presented in
homogeneous coordinates. In the affine model, where
$$
x=\frac{x_1}{x_0}, \qquad y=\frac{y_1}{y_0}
$$
the absolute figure is determined by $w$: $x_0=0$;
$\Omega_1=(0:1:1)$ and
$\Omega_2=(0:1:-1)$.\\
In the pseudo-Euclidean plane the scalar product for two vectors,
e.g. $\mathbf{v_1}=(x_1,y_1)$ and  $\mathbf{v_2}=(x_2,y_2)$, $x_i,
y_i\in \mathbb{R}, i=1,2$ is defined as
\begin{equation}\label{scalar}
\mathbf{v_1} \cdot
\mathbf{v_2}=(x_1,y_1)\cdot(x_2,y_2)=x_1x_2-y_1y_2.
\end{equation}
Hence, the norm of the vector $\mathbf{v}=(x, y)$ is of the form
\begin{equation}\label{norma}
|\mathbf{v}|=\sqrt{\mathbf{v} \cdot \mathbf{v}}=\sqrt{(x, y)\cdot(x,
y)}=\sqrt{x^2-y^2}.
\end{equation}
Since (\ref{norma}) may not always be real, one can distinguish three types of vectors
in the pseudo-Euclidean plane:
\begin{equation}
\begin{tabular}{l}\label{broj}
$1. \quad \textrm{\textit{spacelike vectors} if}\quad \mathbf{v}
\cdot
\mathbf{v}>0;$\\
$2. \quad \textrm{\textit{timelike vectors} if}\quad \mathbf{v}
\cdot
\mathbf{v}<0;$\\
$3. \quad \textrm{\textit{lightlike vectors} (isotropic vectors)
if}\quad \mathbf{v} \cdot
\mathbf{v}=0.$\\
\end{tabular}
\end{equation}
As a consequence there are 3 types of straight lines:
\textit{spacelike lines},
\textit{timelike lines}, \textit{lightlike lines}.\\
Apparently, for two points $T_1=(x_1,
y_1)$ and $T_2=(x_2,y_2)$
\begin{equation}\label{udaljenost}
d(T_1, T_2):=\sqrt{(x_1-x_2)^2-(y_1-y_2)^2}
\end{equation}
defines the distance between them. Comparing (\ref{udaljenost}) and (\ref{norma}), for $\mathbf{v}=\overrightarrow{T_1T_2}$
we have $|\mathbf{v}|=d(T_1,T_2)$. We will use the following notation: $d(T_1,T_2)=|T_1T_2|$.\\
If $0=(0,0)$ is the origin, the vectors $\overrightarrow{OT_1}$ and
$\overrightarrow{OT_2}$, being both spacelike or both timelike, form an angle defined by
\begin{equation}\label{kut}
\cosh
\alpha:=\frac{x_1x_2-y_1y_2}{\sqrt{{x_1}^2-{y_1}^2}\sqrt{{x_2}^2-{y_2}^2}}.
\end{equation}
The transformations that keep the absolute figure invariant and
preserve the above given metric quantities of a scalar product,
distance, angle, are of the form
\begin{equation}\label{gibanje}
\begin{tabular}{l}
$\overline{x}=x\cosh \varphi + y \sinh \varphi + a$\\
$\overline{y}=x \sinh \varphi+y \cosh \varphi+b.$
\end{tabular}
\end{equation}
The transformations (\ref{gibanje}) form a group $B_3$, called the \textit{motion group}. Hence, the group of
pseudo-Euclidean motions consists of \textit{translations} and
\textit{pseudo-Euclidean rotations}, that is
$$
\begin{tabular}{l}
$\overline{x}=x+a$\\
$\overline{y}=y+b$
\end{tabular}
\qquad \textrm{and} \qquad
\begin{tabular}{l}
$\overline{x}=x\cosh \varphi + y \sinh \varphi$\\
$\overline{y}=x \sinh \varphi+y \cosh \varphi$.\\
\end{tabular}
$$
With the geometry of the pseudo-Euclidean plane (also known as \textit{Minkowski plane} and \textit{Lorentzian plane}) one can get
acquainted through, for example, \cite{yag} and \cite{bir}.
\section{Conic equation}
General second-degree equation in two variables can be written in
the form
\begin{equation}\label{konika}
F(x,y)\equiv a_{11}x^2+2a_{12}xy+a_{22}y^2+2a_{01}x+2a_{02}+a_{00}=0
\end{equation}
where $a_{11} \ldots a_{00}\in \mathbb{R}$ and at least one of the
numbers $ a_{11}, a_{12}, a_{22}\neq 0$. All the
solutions of the equation (\ref{konika}) represent the locus of
points in a plane which is called a \textit{conic section} or
simply, a \textit{conic}.\\
Using the matrix notation, we have
\begin{equation}
\begin{tabular}{l}
$F(x,y)\equiv \left[\begin{array}{ccc}1 & x & y \end{array}
\right]\left[\begin{array}{ccc}a_{00} & a_{01} & a_{02}\\ a_{01} & a_{11} & a_{12}\\a_{02} & a_{12} & a_{22}\\
\end{array} \right]\left[\begin{array}{c}1\\ x \\ y
\end{array}\right]=$\\
=$\left[\begin{array}{cc}x & y \end{array} \right] \left[\begin{array}{cc}a_{11} & a_{12}\\ a_{12} & a_{22}\\
\end{array} \right]\left[\begin{array}{c}x\\y \end{array} \right]+2\left[\begin{array}{cc}a_{01} & a_{02} \end{array} \right]\left[\begin{array}{c}x\\y \end{array}
\right]+a_{00}=0$
\end{tabular}
\end{equation}
where
\begin{equation}
A:=\left[\begin{array}{ccc}a_{00} & a_{01} & a_{02}\\ a_{01} & a_{11} & a_{12}\\a_{02} & a_{12} & a_{22}\\
\end{array} \right] \quad \textrm{and} \quad \sigma:= \left[\begin{array}{cc}a_{11} & a_{12}\\ a_{12} & a_{22}\\
\end{array} \right]
\end{equation}
are real, symmetric matrices. In the sequel we will use the
following functions of the coefficients $a_{ij}, i,j=0,1,2$
\begin{equation}\label{invarijante}
\begin{tabular}{l}
$I_1:=a_{11}-a_{22}, \qquad I_2:= det \sigma=\left|\begin{array}{cc}a_{11} & a_{12}\\ a_{12} & a_{22}\\
\end{array} \right|, \qquad I_3:=det A =\left|\begin{array}{ccc}a_{00} & a_{01} & a_{02}\\ a_{01} & a_{11} & a_{12}\\a_{02} & a_{12} & a_{22}\\
\end{array} \right|$\\
$I_4:=\left|\begin{array}{cc}a_{00} & a_{01}\\ a_{01} & a_{11}\\
\end{array} \right|-\left|\begin{array}{cc}a_{00} & a_{02}\\ a_{02} & a_{22}\\
\end{array} \right|, \qquad I_5:=a_{00}.$\\
\end{tabular}
\end{equation}

The aim is to determine the invariants of conics with respect to the
motion group $B_3$ in the pseudo-Euclidean plane. For that purpose,
let's first apply on the conic equation (\ref{konika}) the
''pseudo-Euclidean rotation'' from (\ref{gibanje}) given by:
\begin{equation}\label{rotacija}
\begin{tabular}{l}
$x=\overline{x}\cosh \varphi+\overline{y }\sinh \varphi$\\
$y=\overline{x}\sinh \varphi+\overline{y}\cosh \varphi.$\\
\end{tabular}
\end{equation}
Using matrix notation, (\ref{rotacija}) can be represented as
\begin{equation}\label{rot}
\begin{tabular}{l}
$\left[\begin{array}{c}x\\ y\end{array}
\right]=\left[\begin{array}{cc}\cosh \varphi & \sinh \varphi\\
\sinh \varphi & \cosh \varphi \end{array}
\right]\left[\begin{array}{c}\overline{x}\\ \overline{y}\end{array}
\right],\qquad$ $R:=\left[\begin{array}{cc}\cosh \varphi & \sinh \varphi\\
\sinh \varphi & \cosh \varphi \end{array} \right].$\\
\end{tabular}
\end{equation}
Let's focus on the properties of the matrix $R$ given in
(\ref{rot}):
\begin{itemize}
\item[a)] $\textrm{det}R=\left|\begin{array}{cc}\cosh \varphi & \sinh \varphi\\
\sinh \varphi & \cosh \varphi \end{array} \right|=\cosh^2
\varphi-\sinh^2 \varphi=1$
\item[b)]$R^{-1}=\left[\begin{array}{cc}\cosh \varphi & -\sinh \varphi\\
-\sinh \varphi & \cosh \varphi \end{array} \right]$
\item[c)] $R^T=R$
\item[d)] Denoting columns of $R$ by $\mathbf{v_1}=(\cosh \varphi,\sinh
\varphi)$ and $\mathbf{v_2}=(\sinh \varphi,\cosh \varphi)$ we get
$$
\mathbf{v_1} \cdot \mathbf{v_2}=(\cosh \varphi,\sinh \varphi) \cdot
(\sinh \varphi,\cosh \varphi)=\cosh \varphi \sinh \varphi- \sinh
\varphi \cosh \varphi=0
$$
Computing norms of the vectors $\mathbf{v_1},\mathbf{v_2}$, that is
$$
\begin{tabular}{l}
$|\mathbf{v_1}|=\sqrt{\mathbf{v_1} \cdot \mathbf{v_1}}=\sqrt{(\cosh
\varphi,\sinh \varphi)\cdot (\cosh \varphi,\sinh \varphi)}=$\\
$=\sqrt{\cosh^2\varphi-\sinh^2
\varphi}=\sqrt{1}=1$\\
$ $\\
$|\mathbf{v_2}|=\sqrt{\mathbf{v_2} \cdot
\mathbf{v_2}}=\sqrt{(\sinh \varphi,\cosh \varphi)(\sinh
\varphi,\cosh
\varphi)}=$\\
$=\sqrt{\sinh^2
\varphi-\cosh^2\varphi}=\sqrt{-1}=i,$\\
\end{tabular}
$$
we conclude the columns of $R$ are orthonormal in the
pseudo-Euclidean sense.
\end{itemize}
Because of the aforementioned properties of the matrix $R$, we will say
that $R$ is a \textit{pseudo-orthogonal} matrix.\\
Hence, applying (\ref{rot}) on the conic equation (\ref{konika});
$$
\begin{tabular}{l}
$\left[\begin{array}{cc}\overline{x} & \overline{y }\end{array} \right]\left[\begin{array}{cc}\cosh \varphi & \sinh \varphi\\
\sinh \varphi & \cosh \varphi \end{array} \right] \left[\begin{array}{cc}a_{11} & a_{12}\\ a_{12} & a_{22}\\
\end{array} \right]\left[\begin{array}{cc}\cosh \varphi & \sinh \varphi\\
\sinh \varphi & \cosh \varphi \end{array}
\right]\left[\begin{array}{c}\overline{x}\\\overline{y} \end{array}
\right]+$\\
$+2\left[\begin{array}{cc}a_{01} & a_{02} \end{array}
\right]\left[\begin{array}{cc}\cosh \varphi & \sinh \varphi\\
\sinh \varphi & \cosh \varphi \end{array}
\right]\left[\begin{array}{c}\overline{x}\\\overline{y }\end{array}
\right]+a_{00}=0$
\end{tabular}
$$
one gets
\begin{equation}\label{novak}
F(\overline{x},\overline{y})\equiv \left[\begin{array}{ccc}1 &
\overline{x} & \overline{y}
\end{array}
\right]\left[\begin{array}{ccc}\overline{a_{00}} & \overline{a_{01}} & \overline{a_{02}}\\ \overline{a_{01}} & \overline{a_{11} }& \overline{a_{12}}\\\overline{a_{02}} & \overline{a_{12}} & \overline{a_{22}}\\
\end{array} \right]\left[\begin{array}{c}1\\\overline{ x }\\ \overline{y}
\end{array}\right]=0,
\end{equation}
where
\begin{equation}\label{koeficijenti}
\begin{tabular}{l}
$\overline{a_{11}}=a_{11}{\cosh^2 \varphi}+a_{22}{\sinh^2
\varphi}+2a_{12}\cosh \varphi \sinh \varphi$\\
$\overline{a_{12}}=(a_{11}+a_{22})\cosh \varphi \sinh
\varphi+a_{12}({\cosh^2 \varphi}+{\sinh^2 \varphi})$\\
$\overline{a_{22}}=a_{11}{\sinh^2 \varphi}+a_{22}{\cosh^2
\varphi}+2a_{12}\cosh \varphi \sinh \varphi$\\
$\overline{a_{01}}=a_{01}\cosh \varphi + a_{02}\sinh \varphi$\\
$\overline{a_{02}}=a_{01}\sinh \varphi + a_{02}\cosh \varphi$\\
$\overline{a_{00}}=a_{00}.$\\
\end{tabular}
\end{equation}
This yields $I_1, I_2, I_3, I_4, I_5$ are invariant with
respect to the rotations (\ref{rotacija}).\\
For example,
$$
\begin{tabular}{l}
$\overline{I_3}=\left|\begin{array}{ccc}\overline{a_{00}} & \overline{a_{01}} & \overline{a_{02}}\\ \overline{a_{01}} & \overline{a_{11}} &\overline{ a_{12}}\\\overline{a_{02}} & \overline{a_{12}} & \overline{a_{22}}\\
\end{array} \right|=-\overline{a_{00}} {\overline{a_{12}}}^2 + 2 \overline{a_{01}} \overline{a_{12}} \overline{a_{02}} - \overline{a_{11}} {\overline{a_{02}}}^2 - {\overline{a_{01}}}^2 \overline{a_{22}} + \overline{a_{00}}\overline{
a_{11}} a_{22}=$\\
$-a_{00} {a_{12}}^2 \cosh^4\varphi + 2 a_{01} a_{12} a_{02}
\cosh^4\varphi -
 a_{11} {a_{02}}^2 \cosh^4\varphi - {a_{01}}^2 a_{22} \cosh^4\varphi+
 a_{00} a_{11} a_{22} \cosh^4\varphi +$\\
$+ 2 a_{00} a_{12}^2 \cosh^2\varphi \sinh^2\varphi -
 4 a_{01} a_{12} a_{02} \cosh^2\varphi \sinh^2\varphi+
 2 a_{11} {a_{02}}^2 \cosh^2\varphi \sinh^2\varphi+$\\
$+ 2 a_{01}^2 a_{22} \cosh^2\varphi \sinh^2\varphi -
 2 a_{00} a_{11} a_{22} \cosh^2\varphi \sinh^2\varphi -
 a_{00} {a_{12}}^2 \sinh^4\varphi + 2 a_{01} a_{12} a_{02} \sinh^4\varphi
 -$\\
 $-a_{11} {a_{02}}^2\sinh^4\varphi - {a_{01}}^2 a_{22} \sinh^4\varphi  +
 a_{00} a_{11} a_{22} \sinh^4\varphi=$\\
 $=-a_{00} a_{12}^2 + 2 a_{01} a_{12} a_{02} - a_{11} a_{02}^2 - a_{01}^2 a_{22} + a_{00} a_{11}
 a_{22}=I_3.$\\
 \end{tabular}
$$
The same can be proved for $I_1$, $I_2$, $I_4$ and $I_5$, as well.\\
Taking translations from (\ref{gibanje}) given by
\begin{equation}\label{translacija}
\begin{tabular}{l}
$x=\overline{x}+x_0$\\
$y=\overline{y}+y_0$\\
\end{tabular}
\end{equation}
the equation (\ref{konika}) turns into (\ref{novak}) where
\begin{equation}
\begin{tabular}{l}
$\overline{a_{11}}=a_{11}$\\
$\overline{a_{12}}=a_{12}$\\
$\overline{a_{22}}=a_{22}$\\
$\overline{a_{01}}=a_{11}x_0 + a_{12}y_0+a_{01}$\\
$\overline{a_{02}}=a_{12}{x_0} + a_{22}y_0+a_{02}$\\
$\overline{a_{00}}=a_{11}{x_0}^2+2a_{12}x_0y_0+a_{22}{y_0}^2+2a_{01}x_0+2a_{02}y_0+a_{00}.$\\
\end{tabular}
\end{equation}
It is easy to show that $I_1, I_2, I_3$ are invariants under (\ref{translacija}). One concludes that $I_1$, $I_2$, and $I_3$ are invariants of conics with
respect to the group of motions $B_3$.\\
The observation given above regarding the invariants $I_1, I_2, I_3,
I_4, I_5$ can be found in \cite{rew}. In addition, Reveruk
\cite{rew} defines conics with respect to their relationship to the
absolute figure, relying on the fact that the focus points (foci)
are the points of intersection of the isotropic tangents at the
conic. The paper, however, showed to be incomplete (see Tables 1-5, where the conics added from us are written in italic) and not possible
to cite in further studies of the properties of conics in the
pseudo-Euclidean plane.\\

\section{Diagonalization of the quadratic form}
In the chapters that follows, based on the methods of linear
algebra, we give a complete classification of conic sections, divide
them into families and define types, giving both of them geometrical
meaning.\\
The quadratic form within the equation (\ref{konika}) is a second
degree homogenous polynomial
\begin{equation}\label{kvadrat}
Q(x,y):=a_{11}x^2+2a_{12}xy+a_{22}y^2=\left[\begin{array}{cc}x & y
\end{array} \right] \left[\begin{array}{cc}a_{11} & a_{12}\\a_{12} & a_{22}
\end{array} \right]\left[\begin{array}{c}x\\y
\end{array} \right].
\end{equation}
The question is whether and when it is possible to obtain
$\overline{a_{12}}=0$ using transformations of the group $B_3$. It
can be seen from (\ref{koeficijenti}) that $\overline{a_{12}}=0$
implies
\begin{equation}\label{th}
\begin{tabular}{l}
$(a_{11}+a_{22})\cosh\varphi\sinh\varphi+a_{12}(\cosh^2\varphi+\sinh^2\varphi)=0,$\\
$\displaystyle
\frac{1}{2}(a_{11}+a_{22})\sinh2\varphi+a_{12}\cosh2\varphi=0,$\\
$\textrm{i.e.} \qquad
\tanh2\varphi=-\displaystyle\frac{2a_{12}}{a_{11}+a_{22}}, \quad
a_{11}+a_{22}\neq 0.$
\end{tabular}
\end{equation}
From (\ref{th}) we read:
\begin{itemize}
\item[(i)] $-1<\tanh2\varphi<1, -\infty<2\varphi<\infty $ is
fulfilled when $|a_{11}+a_{22}|>2|a_{12}|$;
\item[(ii)] $\tanh2\varphi=1$, $2\varphi=\infty$ is fulfilled when
$a_{11}+a_{22}=-2a_{12}$,\\
$\tanh2\varphi=-1$, $2\varphi=-\infty$ is fulfilled when
$a_{11}+a_{22}=2a_{12}$,
\item[(iii)] $\tanh2\varphi<-1$ and  $\tanh2\varphi>1$ is
impossible. This follows when $|a_{11}+a_{22}|<2|a_{12}|$.\\
\end{itemize}
So, under the condition (i) one obtains
\begin{equation}\label{kvapotez}
Q(\overline{x}, \overline{y})=\left[\begin{array}{cc}\overline{x} &
\overline{y}
\end{array} \right] \left[\begin{array}{cc}\overline{a_{11}} & 0\\0 &\overline{ a_{22}}
\end{array} \right]\left[\begin{array}{c}\overline{x}\\\overline{y}
\end{array} \right],
\end{equation}
where $\overline{a_{11}}-\overline{a_{22}}=I_1$,
$\overline{a_{11}}\cdot \overline{a_{22}}=I_2$.
\begin{defi}\label{def1}
Let $A:=\left[\begin{array}{cc}a_{11} & a_{12}\\a_{12} & a_{22}
\end{array} \right]$ be any real symmetric matrix. Then the values
$\lambda_1, \lambda_2$,
$$
\lambda_1-\lambda_2=a_{11}-a_{22}, \quad \lambda_1\cdot
\lambda_2=a_{11}a_{22}-{a_{12}}^2
$$
are called pseudo-Euclidean values of the matrix $A$.
\end{defi}
\begin{defi}\label{def2}
We say that the real symmetric $2\times2$ matrix $A$ allows the
pseudo-Euclidean diagonalization if there is a matrix
$$
R=\left[\begin{array}{cc}\cosh\varphi & \sinh\varphi\\ \sinh\varphi
&\cosh\varphi
\end{array} \right]
$$
such that $RAR$ is a diagonal matrix, i.e.
$$
RAR=\left[\begin{array}{cc}\lambda_1& 0\\
0 & \lambda_2
\end{array} \right],
$$
where $\lambda_1, \lambda_2$ are the pseudo-Euclidean values of the
matrix $A$. We say that the matrix $R$ diagonalizes $A$ in a
pseudo-Euclidean way.
\end{defi}
From the results obtained in Section 2 related to the invariants
(\ref{invarijante}) it follows:
\begin{prop}\label{prop1}
The difference  $\lambda_1-\lambda_2$ of the pseudo-Euclidean values
as well as their product $\lambda_1\cdot \lambda_2$ are invariant
with respect to the group $B_3$ of motions in the pseudo-Euclidean
plane.
\end{prop}
Out of (\ref{kvadrat}), (i), (ii), (iii), and (\ref{kvapotez}),
Propositions \ref{prop2} and \ref{prop3} are valid:
\begin{prop}\label{prop2}
Let $A$ be a matrix from Definition \ref{def1}. Then there is a
matrix $R=\left[\begin{array}{cc}\cosh\varphi & \sinh\varphi\\
\sinh\varphi &\cosh\varphi
\end{array} \right]$ with
$\tanh2\varphi=-\displaystyle\frac{2a_{12}}{a_{11}+a_{22}}$ which
under the conditions $a_{11}+a_{22}\neq 0$ and $
|a_{11}+a_{22}|>2|a_{12}|$ diagonalizes $A$ in the pseudo-Euclidean
way.
\end{prop}
\begin{prop}\label{prop3}
It is always possible to reduce the quadratic form (\ref{kvadrat})
by a pseudo-Euclidean motion to the canonical form (\ref{kvapotez})
except for: (ii) and (iii).
\end{prop}
Next we divide the conics in the pseudo-Euclidean plane in four
families according to their geometrical properties. First, we define
the families:
\begin{defi}\label{def3} 1st family conics in the
pseudo-Euclidean plane are conics with no real isotropic directions
while their isotropic points are spacelike
or timelike.\\
2nd family conics are conics having one real isotropic
direction.\\
3rd family conics are conics with two real isotropic points, one being spacelike  and the other being timelike.\\
4th family conics are ones incident with both absolute points.
\end{defi}
Taking in consideration the range of angles in the pseudo-Euclidean
plane \cite{yag}, \cite{shev}, the significance of the conditions (i), (ii), (iii) as well
as that of the equality $a_{11}+a_{22}=0$ is given in the
proposition that follows.
\begin{prop}\label{prop4} Any conic that satisfies the
condition (i) within its equation (\ref{konika}) represents a conic
with no real isotropic directions while their isotropic points are
spacelike or timelike. When one of the conditions (ii) is fulfilled,
(\ref{konika}) represents a conic having one real isotropic
direction. For (iii) (\ref{konika}) represents a conic with two real
isotropic points, one being spacelike and the other being timelike.
Finally, when $a_{11}+a_{22}=0$ is fulfilled, (\ref{konika}) is a
conic incident with both absolute points.
\end{prop}
Let us now discuss the geometrical meaning of the invariants as
follows:
\begin{itemize}
\item[1.] $I_3\neq0$ represents a proper conic while $I_3=0$ represents a degenerate conic.  \vspace{3 mm}
\item[2.] $I_2\neq0$ represents a conic with center and $I_2=0$ a conic without
center. As it is well known 1. and 2. are affine conditions for
conics. \vspace{3 mm}
\item[3.] Conics belonging to the 1st family with $I_1\neq 0$ are conics without real isotropic
directions while those with $I_1=0$ have imaginary isotropic
directions. \vspace{3 mm}
\item[4.] Conics belonging to the 2nd family with $I_1\neq 0$ are conics with one isotropic direction.
If $I_1=0$ is valid the considered conic is a conic with double
isotropic direction. \vspace{3 mm}
\item[5.] Conics belonging to the 4th family with $I_1\neq 0$ are conics with two isotropic
directions. If $I_1=0$ is valid the considered conic is a conic
consisting of an absolute line and one more line.
\end{itemize}
Furthermore, for conics with isotropic points of the same type we have introduced the following notations:
\begin{itemize}
\item[-] \textit{first type conic} is a conic with spacelike isotropic points;
\item[-] \textit{second type conic} is a conic with timelike isotropic points.
\end{itemize}

\section{Pseudo-Euclidean classification of conics}
In anticipation of classifying conics based on their isometric
invariants, we give the pseudo-Euclidean classification based on
families and types of conics in the projective model, in order to
point out the need for our investigation. The projective
representations of conics from \cite{rew} are given in black, while
the ones we have completed Reveruk's classification with are drawn
in red color (see figures 1, 2, and 3).
\begin{figure}\centering
\includegraphics[height=8.5cm]{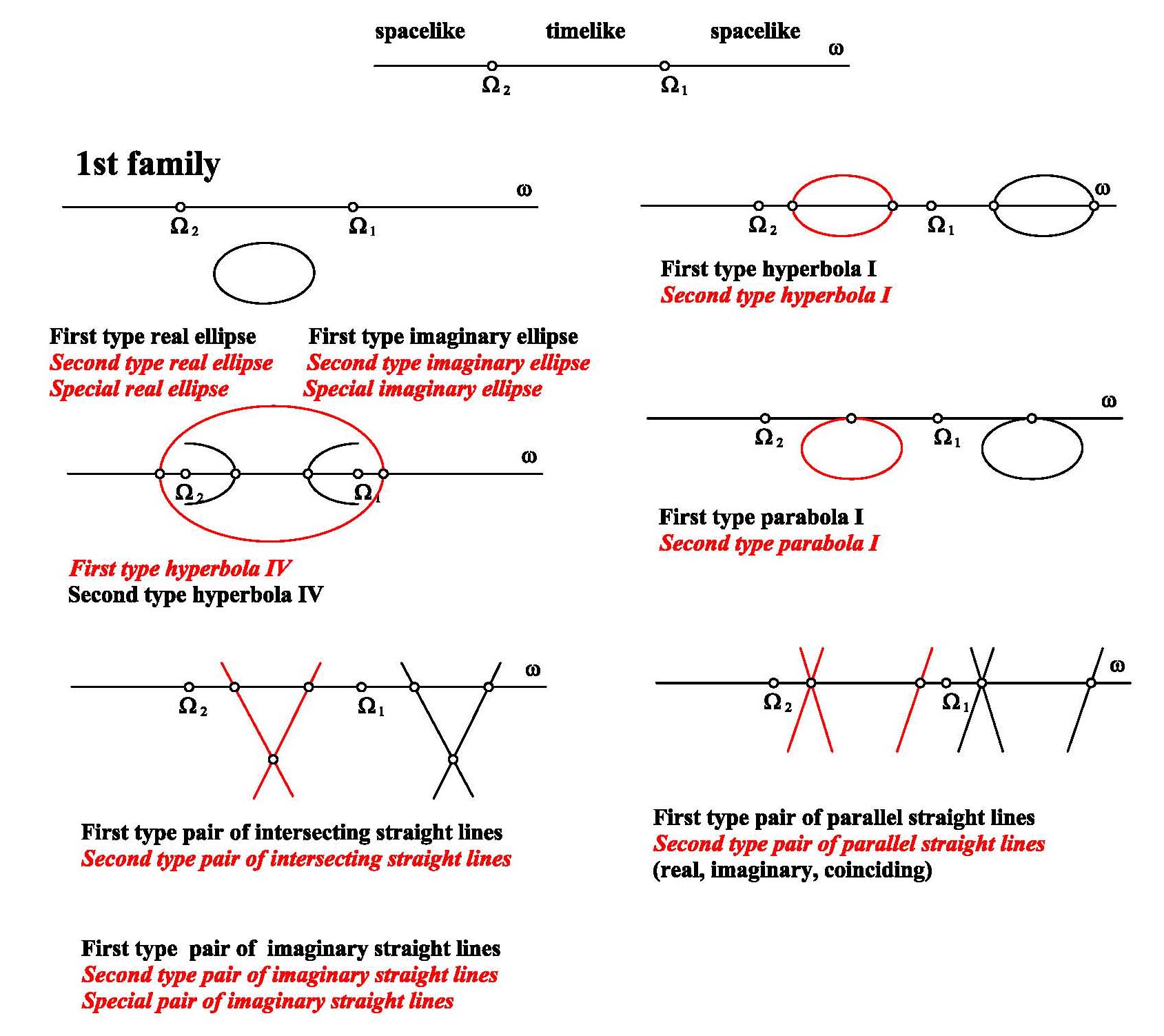}\\
\caption{1st family conics}
\end{figure}
\begin{figure}\centering
\includegraphics[height=8.5cm]{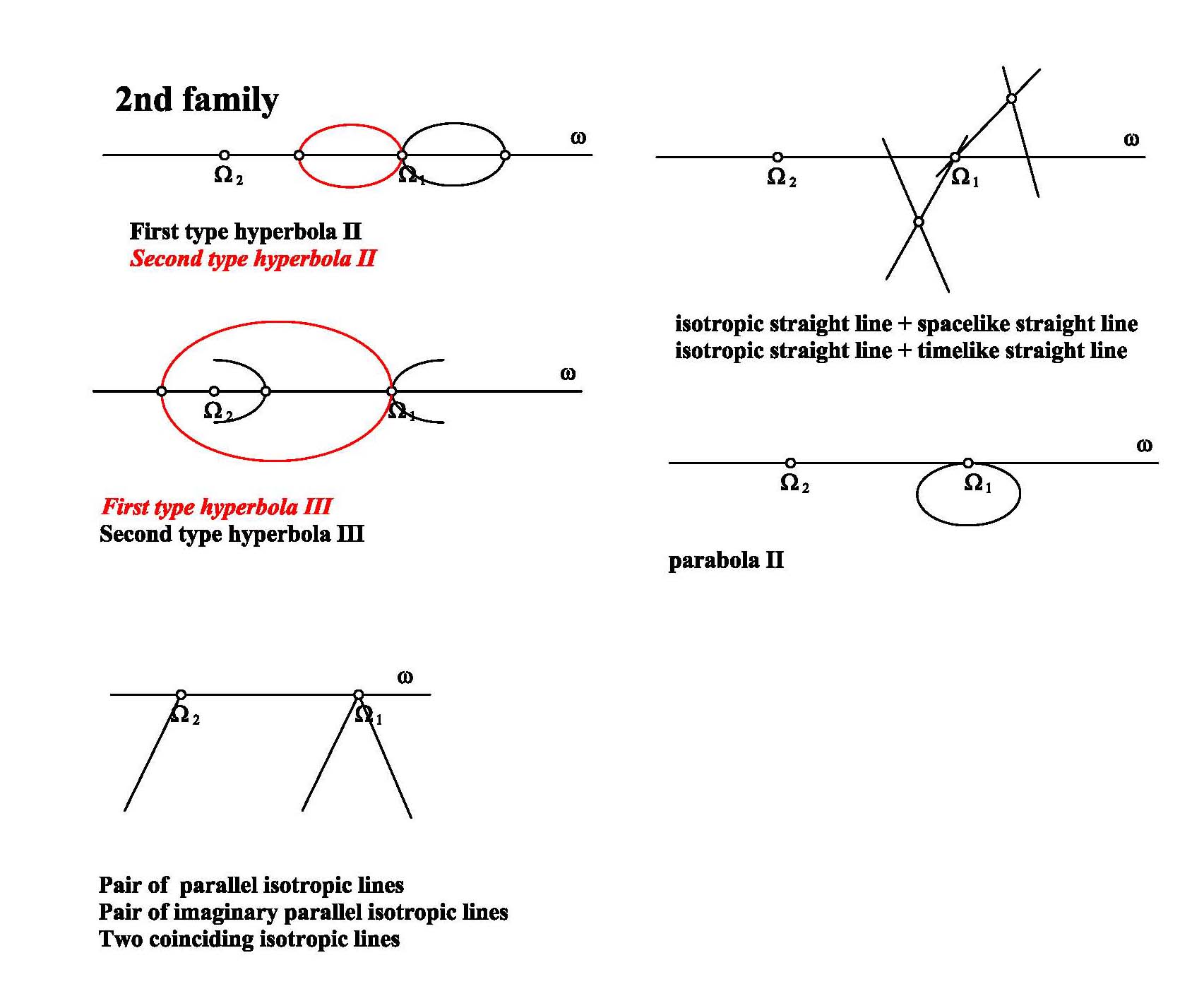}\\
\caption{2nd family conics}
\end{figure}
\begin{figure}\centering
\includegraphics[height=8.5cm]{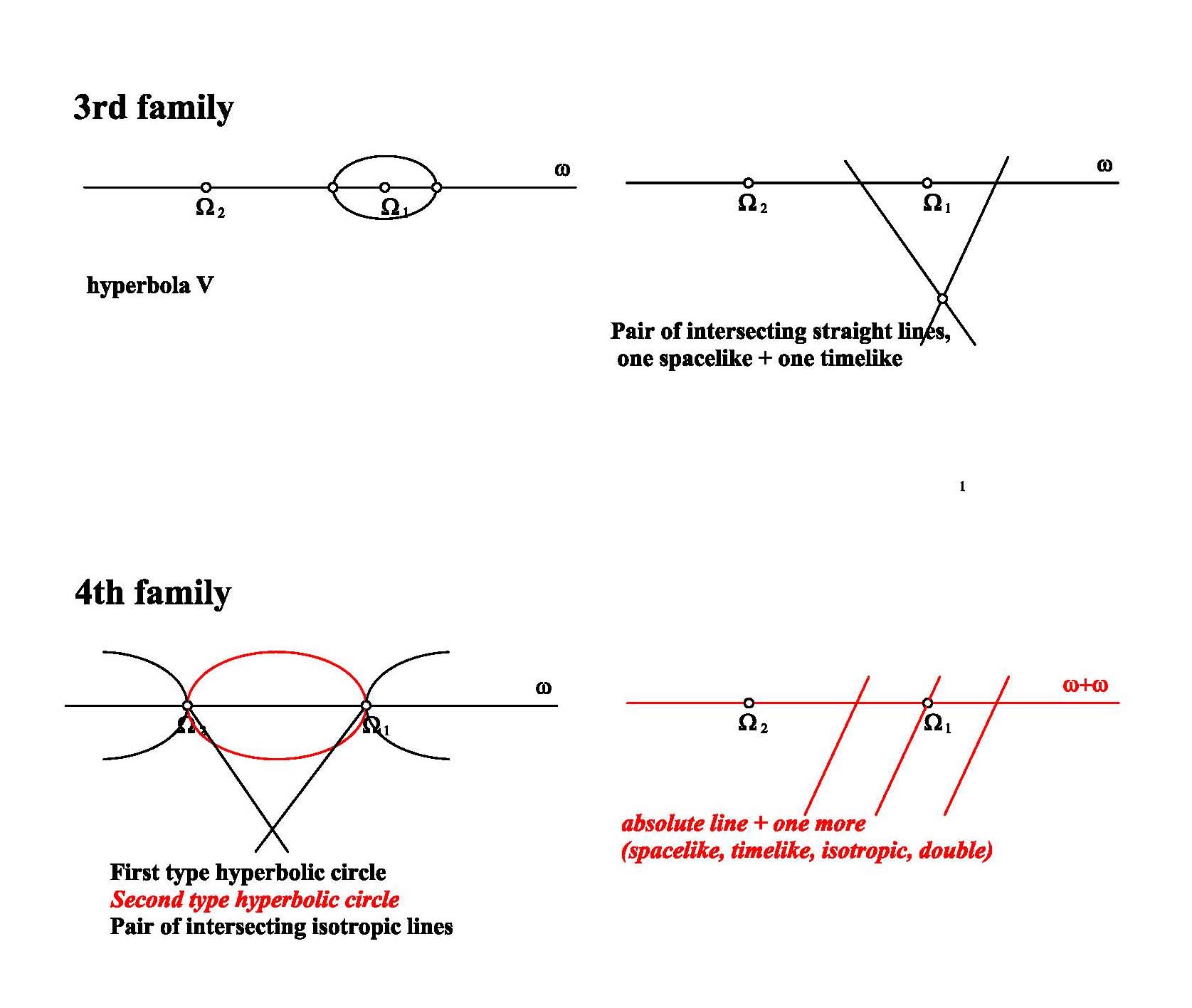}\\
\caption{3rd and 4th family conics}
\end{figure}
\subsection{1st family conics}
Let's assume that it is possible to reduce the quadratic form in the
conic equation (\ref{konika}) to the canonical form
(\ref{kvapotez}). This implies according to Propositions \ref{prop2}
and \ref{prop3} that $|a_{11}+a_{22}|>2|a_{12}|$, and that it is
possible to write down the conic equation (\ref{konika}) in the form
\begin{equation}\label{konpotez}
F(\overline{x},\overline{y})\equiv
\overline{a_{11}}\overline{x}^2+\overline{a_{22}}\overline{y}^2+2\overline{a_{01}}\overline{x}+2\overline{a_{02}}\overline{y}+\overline{a_{00}}=0.
\end{equation}
Let's consider conics with center ($I_2\neq0$).\\
After a
translation of the coordinate system in $\overline{x}$-and
$\overline{y}$-direction we have
\begin{equation}\label{prva}
F(\overline{x},\overline{y})\equiv
\overline{a_{11}}\overline{x}^2+\overline{a_{22}}\overline{y}^2+\overline{\overline{a_{00}}}=0.
\end{equation}
One computes
\begin{equation}\label{invarpotez}
I_1=\overline{a_{11}}-\overline{a_{22}},\quad
I_2=\overline{a_{11}}\cdot \overline{a_{22}}, \quad
I_3=\overline{a_{11}} \cdot \overline{a_{22}} \cdot
\overline{\overline{a_{00}}} \Rightarrow
\overline{\overline{a_{00}}}=\frac{I_3}{I_2}.
\end{equation}
Let's introduce:
\begin{equation}\label{poluosi}
a:=\sqrt{\left|\frac{\overline{\overline{a_{00}}}}{\overline{a_{11}}}\right|},
\qquad
b:=\sqrt{\left|\frac{\overline{\overline{a_{00}}}}{\overline{a_{22}}}\right|}.
\end{equation}
The values $a$ and $b$ shall be called \textit{pseudo-Euclidean
semiaxes}.\\
In the table that follows we give the possibilities for the conic
sections with equation (\ref{prva}) depending on the signs of the
coefficients. The italic cases are those we added to Reveruk's
classification.\vspace{75mm}\\
Table 1: \vspace{5mm}\\
\begin{tabular}{|l|l|l|l|}
\hline
\tiny{ $\overline{a_{11}}\overline{a_{22}} \overline{a_{00}}$} & \tiny{canonical form} & \tiny{conic}\\
\hline
 \begin{tabular}{l}
 \tiny{$+ + +$}\\
\tiny{$- - -$}
 \end{tabular} &
 \tiny{$\displaystyle \frac{x^2}{a^2}+ \frac{y^2}{b^2}=-1$} &
 \begin{tabular}{l}
\tiny{first type imaginary ellipse ($a>b$)}\\
\tiny{\textit{second type imaginary ellipse} ($a<b$)}\\
\tiny{\textit{special imaginary ellipse} ($a=b$)}
\end{tabular}\\
\hline
 \begin{tabular}{l}
\tiny{$+ + - $}\\
 \tiny{$- - + $}
\end{tabular}
& \tiny{$\displaystyle \frac{x^2}{a^2}+ \frac{y^2}{b^2}=1$} &
\begin{tabular}{l}
\tiny{first type real ellipse ($a>b$)}\\
\tiny{\textit{second type real ellipse} ($a<b$)}\\
\tiny{\textit{special real ellipse} ($a=b$)}
\end{tabular}\\
\hline
\begin{tabular}{l}
\tiny{$+ - -$}\\
\tiny{$- + +$}
 \end{tabular} &
\tiny{$\displaystyle \frac{x^2}{a^2}- \frac{y^2}{b^2}=1$}&
\begin{tabular}{l}
\tiny{first type hyperbola I ($a>b$)}\\
\tiny{second type hyperbola IV ($a<b$)}
 \end{tabular}\\
\hline
\begin{tabular}{l}
\tiny{$- + -$}\\
\tiny{$+ - +$}
 \end{tabular} &
\tiny{$\displaystyle -\frac{x^2}{a^2}+ \frac{y^2}{b^2}=1$} &
\begin{tabular}{l}
\tiny{\textit{second type hyperbola I} ($a<b$)}\\
\tiny{\textit{first type hyperbola IV }($a>b$)}
\end{tabular}\\
\hline
\begin{tabular}{l}
\tiny{$+ + 0$}\\
\tiny{$- - 0$}
 \end{tabular} &
\tiny{$a_{11}x^2+a_{22}y^2=0$} &
\begin{tabular}{l}
\tiny{first type pair of imaginary straight lines ($|a_{11}|<|a_{22}|$)}\\
\tiny{\textit{second type pair of imaginary straight lines }($|a_{11}|>|a_{22}|$)}\\
\tiny{\textit{special pair of imaginary straight lines}
($|a_{11}|=|a_{22}|$)}
\end{tabular}\\
\hline
\begin{tabular}{l}
\tiny{$+ - 0$}\\
\tiny{$- + 0$}
 \end{tabular} &
\tiny{$a_{11}x^2+a_{22}y^2=0$} &
\begin{tabular}{l}
\tiny{first type pair of intersecting straight lines ($|a_{11}|<|a_{22}|$)}\\
\tiny{\textit{second type pair of intersecting straight lines}
($|a_{11}|>|a_{22}|$)}
\end{tabular}\\
\hline
\end{tabular}\vspace{5 mm}\\

The question that naturally arises is why the curves with the
canonical equations given in Table 1. in the pseudo - Euclidean
plane are called as it is given in the same table and what is the
connection between the signs of the coefficients and the conditions
based on the invariants (\ref{invarijante}). We answer by
demonstrating on the case of hyperbola I the procedure conducted for
all the curves from this family.
\subsubsection{First and second type hyperbola I}
\begin{defi}\label{hiperbola} The locus of points in
the pseudo-Euclidean plane for which the difference of their
distances from two different fixed points (foci) in this plane is constant
will be called hyperbola I.
\end{defi}
We distinguish two cases: first and second type hyperbola I.\\
Let $F_1=(c, 0), F_2=(-c, 0), F_3=(0,c), F_4=(0, -c), c\neq 0$ be the given
points. For any point $M=(x,y)$ for which $\overrightarrow{F_1M}$
and $\overrightarrow{F_2M}$ are spacelike vectors, according to
definition \ref{hiperbola}
\begin{equation}\label{novohip}
|F_1M| - |F_2M|=2a, \quad a\in\mathbb{R}, a\neq 0,
\end{equation}
\begin{equation}\label{hipdrugo}
\textrm{i.e.} \qquad \sqrt{(x-c)^2-y^2}-\sqrt{(x+c)^2-y^2}=2a.
\end{equation}
After computing (\ref{hipdrugo}) we get
\begin{equation}\label{hip}
\frac{x^2}{a^2}-\frac{y^2}{b^2}=1, \quad b^2=a^2-c^2, a>b, a>c.
\end{equation}
Out of which, according to the affine classification of the second
order curves, because of $a_{11}a_{22}-a_{12}^2=-a^2b^2<0,
a^4b^4\neq 0$, we conclude that the considered conic is a hyperbola.
The symbol I denotes that the foci are real points, i. e., the points $(0:1:1)$
and $(0:1:-1)$ are lying outside the
hyperbola. It is easy to check that the isotropic points are spacelike points, being property of a first type conic and achieved when $a>b$.\\
Equation (\ref{hip}) can be obtained in much the same way carrying
out a calculation for the points $F_3$ and $F_4$,
$\overrightarrow{F_3M}$ and $\overrightarrow{F_4M}$ being
again spacelike vectors, i. e. $|F_3M| - |F_4M|=2b, \quad b\in \mathbb{R}, b\neq 0$.\\
It is easy to show the opposite direction of the above statement as well, i.e., for any
point $M(x, y)$ whose coordinates fulfill the equation (\ref{hip})
the equality $|F_1M| - |F_2M|=2a$ is valid, i.e. the point $M$ is
incident
to the hyperbola I.\vspace {4 mm}\\
Let's presume next $\overrightarrow{F_1M}$ and
$\overrightarrow{F_2M}$ are timelike vectors,
\begin{equation}\label{imhip}
|F_1M| - |F_2M|=2ai, \quad a\in\mathbb{R}, a\neq 0,
\end{equation}
\begin{equation}\label{kika}
\textrm{i.e.} \quad \sqrt{(x-c)^2-y^2}-\sqrt{(x+c)^2-y^2}=2ai
\end{equation}
From (\ref{kika}) we get
\begin{equation}
-\frac{x^2}{a^2}+\frac{y^2}{b^2}=1, \quad a^2+c^2=b^2, b>a,
\end{equation}
being a hyperbola I, of the second type.\vspace {4 mm}\\
The connection between the signs of the coefficients in the
canonical forms of the discussed conics and the (meeting) conditions
based on the invariants $I_1, I_2, I_3$ is given next:\\
For first type hyperbola I the signs of the coefficients
$\overline{a_{11}}$, $\overline{a_{22}}$,
$\overline{\overline{a_{00}}}$ are $+$, $-$, $-$ or $-$, $+$, $+$,
respectively, and $a>b$. This results in
$$
\begin{tabular}{l}
$I_2<0 \quad \wedge \quad ((I_1>0 \wedge I_3>0)\vee (I_1<0 \wedge
I_3<0))\quad \wedge \quad |\overline{a_{11}}|<|\overline{a_{22}}|$\\
$ $\\
$\textrm{i.e.} \quad I_2<0 \wedge I_1I_3>0 \wedge
|\overline{a_{11}}|<|\overline{a_{22}}|.$
\end{tabular}
$$
The opposite direction holds as well.\\
For second type hyperbola I the signs for $\overline{a_{11}}$,
$\overline{a_{22}}$, $\overline{\overline{a_{00}}}$ are $-$, $+$,
$-$ or $+$, $-$, $+$, and $a<b$. This results in
$$
\begin{tabular}{l}
$I_2<0 \quad \wedge \quad ((I_1>0 \wedge I_3<0)\vee (I_1<0 \wedge
I_3>0))\quad \wedge \quad |\overline{a_{11}}|>|\overline{a_{22}}|$\\
$ $\\
$\textrm{i.e.} \quad I_2<0 \wedge I_1I_3<0 \wedge
|\overline{a_{11}}|>|\overline{a_{22}}|.$
\end{tabular}
$$
Conics of the 1st family with $I_2=0$ may be considered in a similar
way. They are included in Table 5. If it is deemed necessary, those
cases can be discuss as well.\\
We conclude subsection 4.1 with the following proposition:
\begin{prop}\label{profirst}
In the pseudo-Euclidean plane there are 23 (12 proper + 11
degenerate) different types of conic sections of the 1st family to
distinguish with respect to the group $B_3$ of motions (see Table
5).
\end{prop}
\subsection{2nd family conics}
Let's assume furtheron that it is not possible to diagonalize the
quadratic form in the conic equation. Then according to Proposition
\ref{prop3} we have to distinguish (ii) $|a_{11}+a_{22}|=2|a_{12}|$
and (iii) $|a_{11}+a_{22}|<2|a_{12}|$, that is 2nd and 3rd
family conics. \vspace{4 mm}\\
The conditions $|a_{11}+a_{22}|=2|a_{12}|$ and $a_{12}\neq 0$ imply
$a_{11}+a_{22}\neq 0$. The conic equation is of the initial form (\ref{konika}). \vspace {3mm}\\
Let's consider conics with center ($I_2\neq 0$).\\
After a translation of a coordinate system in $\overline{x}$- and
$\overline{y}$- direction we have
\begin{equation}\label{novo1}
F(\overline{x},\overline{y})\equiv
a_{11}\overline{x}^2+2a_{12}\overline{x}\overline{y}+a_{22}\overline{y}^2+\overline{a_{00}}=0.
\end{equation}
One computes
$$
\overline{a_{00}}=\frac{I_3}{I_2}.
$$
The possibilities for the conic sections with equation (\ref{novo1}) are: \vspace{5 mm}\\
Table 2: \vspace{5 mm}\\
\begin{tabular}{|l|l|l|l|}
\hline
\tiny{ $a_{11}a_{12}a_{22}\overline{a_{00}}$} & \tiny{canonical form} & \tiny{conic}\\
 \hline
 \begin{tabular}{l}
\tiny{$+ + - -$}\\
 \tiny{$- - + +$} \vspace{1mm}\\
 \tiny{$+ - - -$}\\
 \tiny{$- + + +$}
\end{tabular}
&
\begin{tabular}{l}
\tiny{$x^2(a^2-c^2) + 2xyc^2-y^2(a^2+c^2)-a^4=0$}\\
$ $\\
\tiny{$x^2(a^2-c^2) - 2xyc^2-y^2(a^2+c^2)-a^4=0$}
\end{tabular} &
\tiny{first type hyperbola II}\\
\hline
 \begin{tabular}{l}
 \tiny{$+ - - +$}\\
\tiny{$- + + -$} \vspace{1mm}\\
 \tiny{$+ + - +$}\\
 \tiny{$- - + -$}
 \end{tabular} &
 \begin{tabular}{l}
 \tiny{$x^2(a^2+c^2) - 2xyc^2-y^2(a^2-c^2)+a^4=0$}\\
 $$\\
 \tiny{$x^2(a^2+c^2)+ 2xyc^2-y^2(a^2-c^2)+a^4=0$}
 \end{tabular} &
\tiny{\textit{second type hyperbola II}}\\
\hline
\begin{tabular}{l}
\tiny{$+ + - +$}\\
\tiny{$- - + -$} \vspace{1mm}\\
\tiny{$+ - - +$}\\
\tiny{$- + + -$}
 \end{tabular} &
 \begin{tabular}{l}
\tiny{$x^2(a^2-c^2) + 2xyc^2-y^2(a^2+c^2)+a^4=0$}\\
$$\\
\tiny{$x^2(a^2-c^2) - 2xyc^2-y^2(a^2+c^2)+a^4=0$}
\end{tabular}
& \tiny{\textit{first type hyperbola III}}\\
\hline
\begin{tabular}{l}
\tiny{$+ - - -$}\\
\tiny{$- + + +$}\vspace{1mm}\\
\tiny{$+ + - -$}\\
\tiny{$- - + +$}
 \end{tabular} &
 \begin{tabular}{l}
 \tiny{$x^2(a^2+c^2) - 2xyc^2-y^2(a^2-c^2)-a^4=0$}\\
 $$\\
 \tiny{$x^2(a^2+c^2)+ 2xyc^2-y^2(a^2-c^2)-a^4=0$}
 \end{tabular} & \tiny{second type hyperbola III}\\
\hline
\begin{tabular}{l}
\tiny{ $+ + - 0$}\\
\tiny{ $- - + 0$}
 \end{tabular} &
 \tiny{$x^2(a^2-c^2) + 2xyc^2-y^2(a^2+c^2)=0$} & \tiny{pair of lines, one isotropic, one spacelike}\\
\hline
\begin{tabular}{l}
\tiny{$+ - - 0$}\\
\tiny{$- + + 0$}
 \end{tabular} &
\tiny{$x^2(a^2+c^2) - 2xyc^2-y^2(a^2-c^2)=0$}& \tiny{pair of lines, one isotropic, one timelike}\\
\hline
\end{tabular}\vspace{5 mm}\\
We point out that Reveruk makes difference by name but not by the
invariants between the degenerate conics from Table 2, as well as
between hyperbolas II and III from the same table.\\
\subsubsection{First and second type hyperbola II}
Let us next turn our attention to, for example, hyperbolas II. We
will demonstrate how their names has been derived from their
canonical equations. In addition we provide a link between the signs
of the coefficients within their canonical equations and the
conditions based on the invariants (\ref{invarijante}) for a conic
to represent first, i. e. second type hyperbola II.
\begin{defi}
The locus of points in the pseudo-Euclidean plane for which the
difference from two fixed points (foci) lying on one of the
isotropic lines is constant is called hyperbola II.
\end{defi}
We distinguish 2 cases: first and second type hyperbola II.\\
Let $F_1=(c, c), F_2=(-c,-c)$ be the given points. For any point
$M=(x,y)$  for which $\overrightarrow{F_1M}$ and
$\overrightarrow{F_2M}$ are spacelike vectors,
\begin{equation}\label{hip1}
|F_1M| - |F_2M|=2a, \quad a\in\mathbb{R}, a\neq 0,
\end{equation}
\begin{equation}\label{hip2}
\textrm{i.e.} \qquad
\sqrt{(x-c)^2-(y-c)^2}-\sqrt{(x+c)^2-(y+c)^2}=2a.
\end{equation}
After computing (\ref{hip2}) we get
\begin{equation}\label{hip3}
x^2(a^2-c^2) + 2xyc^2-y^2(a^2+c^2)-a^4=0, \quad a>c.
\end{equation}
Out of (\ref{hip3}), according to the affine classification of the
second order curves
$a_{11}a_{22}-{a_{12}}^2=-(a^2-c^2)(a^2+c^2)-c^4=-a^4<0$,
$(a_{11}a_{22}-{a_{12}}^2)a_{00}=a^8\neq 0$; it is a matter of a
hyperbola \cite{anton}. Further, $\Omega_1=(0:1:1)$ is lying on while
$\Omega_2=(0:1:-1)$ is lying outside the hyperbola, being properties
of II. As the second isotropic point of the curve belongs to the
spacelike area, it is a matter of a first type curve.\\
Carrying out a calculation for the points $F_1=(-c, c), F_2=(c,
-c)$, lying on the isotropic line $x+y=0$, one gets
\begin{equation}\label{hip4}
x^2(a^2-c^2)-2xyc^2-y^2(a^2+c^2)-a^4=0, \quad a>c,
\end{equation}
being again first type hyperbola II. \vspace {4 mm}\\
Presuming that for
$F_1=(c,c)$, $F_2=(-c, -c)$ and $M=(x,y)$ $\overrightarrow{F_1M}$
and $\overrightarrow{F_2M}$ are timelike vectors, we start from
\begin{equation}
|F_1M| - |F_2M|=2ai
\end{equation}
which leads to
\begin{equation}\label{nesto2}
x^2(a^2+c^2)-2c^2xy-y^2(a^2-c^2)+a^4=0, \quad a>c.
\end{equation}
The equation (\ref{nesto2}) represents a \textit{second type hyperbola II}.\\
Repeating the
calculation for $F_1=(-c,c)$, $F_2(c,-c)$ we get
\begin{equation}
x^2(a^2+c^2)+2c^2xy-y^2(a^2-c^2)+a^4=0,
\end{equation}
being again a \textit{second type hyperbola II}.\\
Same as in the case of hyperbolas I, the opposite direction holds as well.\\
If we discuss the signs of the coefficients for first type hyperbola
II we get: there are two possibilities for the signs of the
coefficients $a_{11}$, $a_{12}$, $a_{22}$,
$\overline{ a_{00}}$\\
$$
(\begin{tabular}{l}
$++--$\\
$--++$
\end{tabular}
\quad \textrm{or} \quad
\begin{tabular}{l}
$+---$\\
$-+++$
\end{tabular})
\quad \textrm{and} \quad |a_{11}|<|a_{22}|.
$$
Both combinations of signs yield
$$
I_2<0\quad \wedge \quad ((I_1>0\wedge I_3>0) \vee (I_1<0 \wedge
I_3<0))
$$
which result in
$$
I_2<0, \quad  I_1I_3>0, \quad  |a_{11}|<|a_{22}|.
$$
For second type hyperbola II we start from
$$
(\begin{tabular}{l}
$+--+$\\
$-++-$
\end{tabular}
\quad \textrm{or} \quad
\begin{tabular}{l}
$++-+$\\
$--+-$
\end{tabular})
\quad \textrm{and} \quad |a_{11}|>|a_{22}|,
$$
which leads to
$$
I_2<0, \quad  I_1I_3>0, \quad  |a_{11}|>|a_{22}|.
$$
The opposite direction is valid in both cases. In a very similar way
conics of the 2nd family with $I_2=0$ are considered. We conclude
the analysis within this family with
\begin{prop}
In the pseudo-Euclidean plane there are 10 (5 proper + 5 degenerate)
different types of conic sections of the 2nd family to distinguish
with respect to the group $B_3$ of motions (see Table 5.).
\end{prop}
\subsection{3rd family conics}
Conic sections of the 3rd family are those with two real isotropic
points, one being spacelike and the other being timelike. Due to
this property conics has to be with a center, of hyperbolic type,
which is provided by $I_2<0$. Apart from that according to
Proposition \ref{prop4} the condition $|a_{11}+a_{22}|<2|a_{12}|$
has to be fulfilled within equation (\ref{konika}).\\
After a translation
$$
\overline{x}=x-\frac{a_{12}a_{02}-a_{22}a_{01}}{a_{11}a_{22}-{a_{12}}^2},
\qquad \overline{y}=
y-\frac{a_{12}a_{01}-a_{11}a_{02}}{a_{11}a_{22}-{a_{12}}^2}
$$
of the coordinate system in $x-$ and $y-$direction, obtained from
$\displaystyle \frac{\partial F}{\partial x}=0$ and $\displaystyle
\frac{\partial F}{\partial y}=0$, for the conic equation
(\ref{konika}) we have
\begin{equation}\label{kon1}
F(\overline{x}, \overline{y})\equiv
a_{11}\overline{x}^2+2a_{12}\overline{x}\overline{y}+a_{22}\overline{y}^2+\overline{a_{00}},
\end{equation}
where $\overline{a_{00}}= \displaystyle \frac{I_3}{I_2}$.\\
The possibilities for the conic sections with the equation
(\ref{kon1}), according to \cite{rew} are: \vspace{5mm}\\
Table 3: \vspace{5mm}\\

\begin{tabular}{|l|l|l|l|}
\hline
\tiny{ $a_{11} a_{12} a_{22} \overline{a_{00}}$} & \tiny{canonical form} & \tiny{conic}\\
\hline
 \begin{tabular}{l}
 \tiny{$+ - + -$}\\
\tiny{$- + - +$}
 \end{tabular} &
 \tiny{$(a^2-c^2)x^2-2(a^2+c^2)xy+(a^2-c^2)y^2-a^4=0$}
 &
 \tiny{hyperbola V}\\
\hline
 \begin{tabular}{l}
\tiny{$+ - + 0$}\\
 \tiny{$- + - 0$}
\end{tabular}
& \tiny{$(a^2-c^2)x^2-2(a^2+c^2)xy+(a^2-c^2)y^2=0$}
 &
\begin{tabular}{l}
\tiny{pair of lines, one spacelike }\\
\tiny{and one timelike}
\end{tabular}\\
\hline
\end{tabular}\vspace{5 mm}\\
As we didn't have to interfere in Reveruk's classification
concerning conics of the 3rd family, for details on obtaining the
canonical forms in Table 3 one can consult \cite{rew}.\\
However, we note that in this case the foci of a hyperbola are
complex conjugate, and in order to comply the canonical form of a
hyperbola with those of the hyperbolas of the 1st and 2nd family for
the asymptotes were selected straight lines of the form
\begin{equation}
x(a-c)-y(a+c)=0, \qquad x(a+c)-y(a-c)=0.
\end{equation}
For the conditions based on the invariants (\ref{invarijante}) to
represent conics of this family see Table 5.
\begin{prop}
In the pseudo-Euclidean plane there are 2 (1 proper + 1 degenerated)
different types of conic sections of the 3rd family to distinguish
with respect to the group $B_3$ of motions (see Table 5.).
\end{prop}
\subsection{4th family conics}
For a complete classification of conic sections in the
pseudo-Euclidean plane it is necessary to take into account the
conic sections incident with both absolute points. According to
Definition \ref{def3}, such curves belong to the 4th family. On the
other hand, according to Proposition \ref{prop4}, the condition
$a_{11}+a_{22}=0$ has to be fulfilled within the conic equation
(\ref{konika}).\\
The conic section equation in homogeneous coordinates
$(x_0:x_1:x_2)$ is of the form
\begin{equation}\label{conich}
F(x_0, x_1, x_2)\equiv
a_{11}{x_1}^2+2a_{12}x_1x_2+a_{22}{x_2}^2+2a_{01}x_1x_0+2a_{02}x_2x_0+a_{00}{x_0}^2=0.
\end{equation}
From the requirement that the conic with equation (\ref{conich}) is incident with the
absolute points $\Omega_1=(0:1:1)$ and $\Omega_2=(0:1:-1)$ is easy
to show that, apart from $a_{11}+a_{22}=0, a_{12}=0$ holds as well.
The equation (\ref{konika}) now turns into
\begin{equation}\label{kon2}
F(x, y)\equiv a_{11}x_{1}^2+a_{22}y^2+2a_{01}x+2a_{02}y+a_{00}=0
\end{equation}
Presuming that $I_2\neq 0$, both linear terms can be eliminated by a
translation in direction of the $x-$ and $y-$ axes, which gives us
\begin{equation}\label{kon3}
F(\overline{x},\overline{y})\equiv
a_{11}\overline{x}^2+a_{22}\overline{y}^2+\overline{a_{00}}=0.
\end{equation}
One computes
\begin{equation}\label{inv4}
I_1=a_{11}-a_{22}, \quad I_2=a_{11}a{22}, \quad
I_3=a_{11}a_{22}\overline{a_{00}}\Rightarrow
\overline{a_{00}}=\frac{I_3}{I_2}.
\end{equation}
The possibilities for the conic sections with equation (\ref{kon3})
are:\vspace{4 mm}\\
Table 4: \vspace{4 mm}\\
\begin{tabular}{|l|l|l|l|}
\hline
\tiny{ $a_{11} a_{22} \overline{a_{00}}$} & \tiny{canonical form} & \tiny{conic}\\
\hline
 \begin{tabular}{l}
 \tiny{$+ - + $}\\
\tiny{$- + -$}
 \end{tabular} &
 \tiny{$x^2-y^2+a^2=0$} &
 \tiny{\textit{second type hyperbolic circle}}\\
\hline
 \begin{tabular}{l}
\tiny{$+ - -$}\\
 \tiny{$- + +$}
\end{tabular}
& \tiny{$x^2-y^2-a^2=0$} &
\tiny{first type hyperbolic circle}\\
\hline
 \begin{tabular}{l}
\tiny{$+ - 0$}\\
 \tiny{$- + 0$}
\end{tabular}
& \tiny{$x^2-y^2=0$} &
\tiny{pair of isotropic lines}\\
\hline
\end{tabular}\vspace{5 mm}\\
Links among the canonical equations and the corresponding names of
conics from Table 4. are obvious. For the conditions based on the
invariants (\ref{invarijante}) to represent those conics see Table
5.\\
\begin{figure} \centering
\includegraphics[height=16cm]{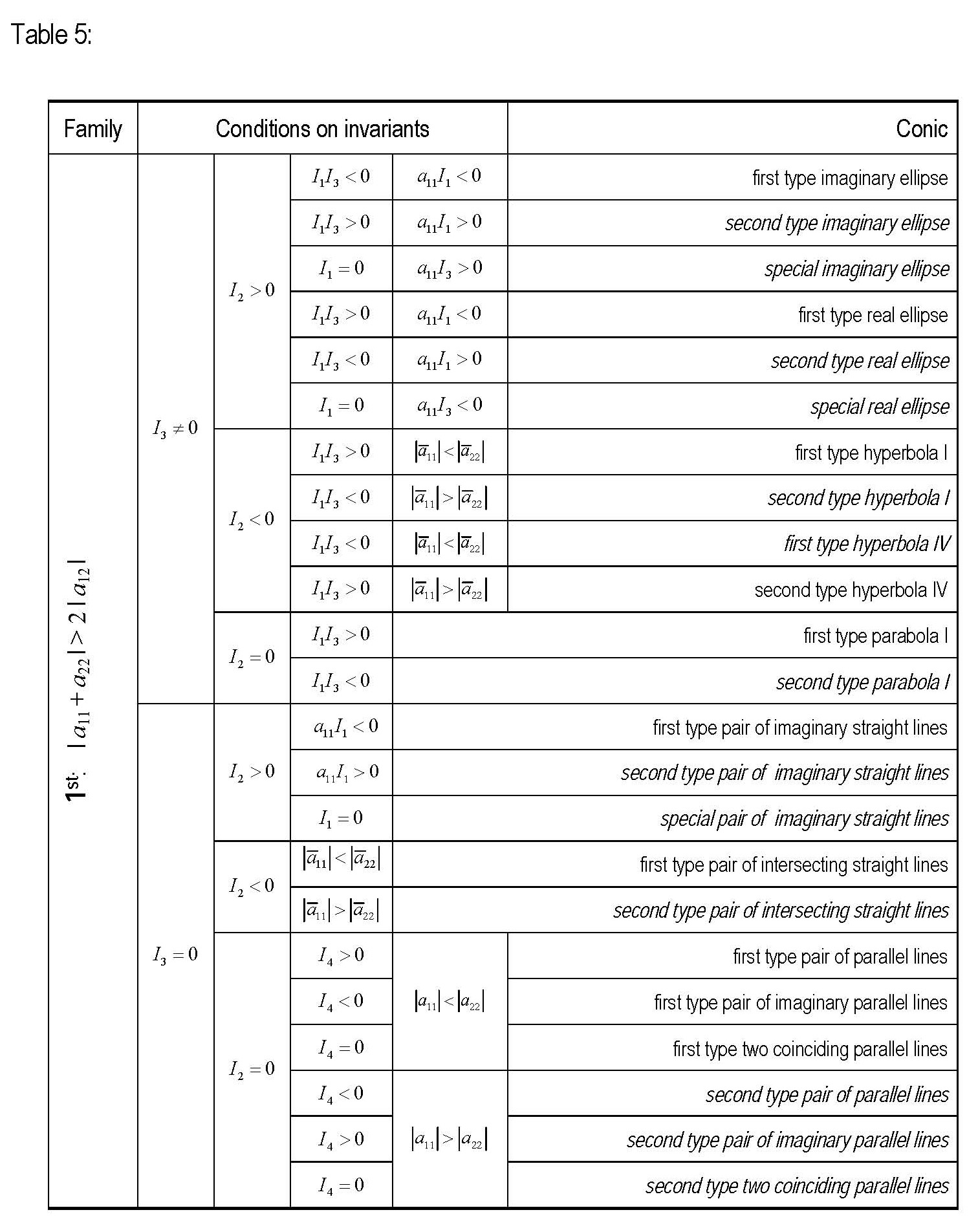}\\
\end{figure}
\begin{figure} \centering
\includegraphics[height=16cm]{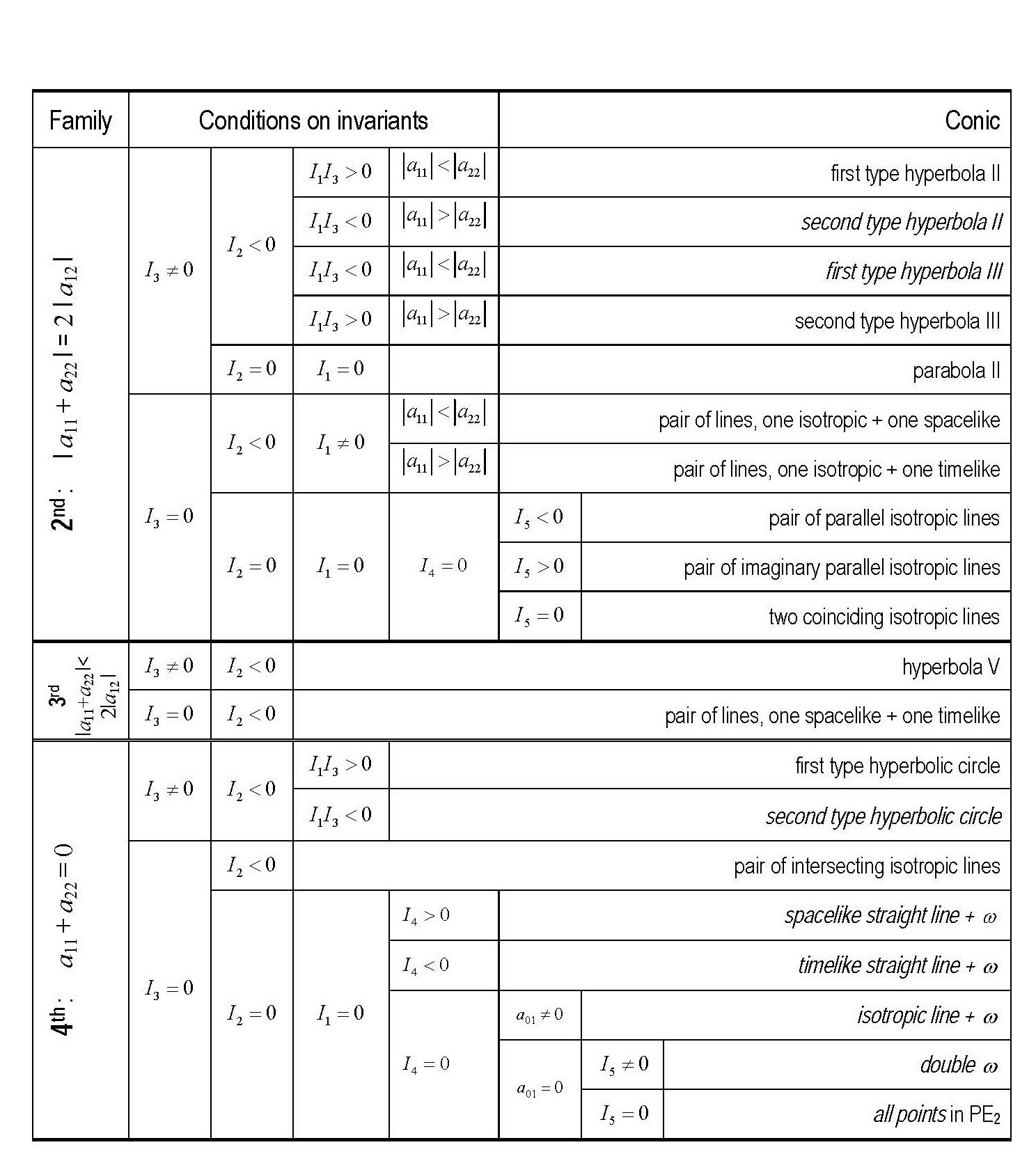}\\
\end{figure}

We continue our study by analyzing conics consisting of two straight
lines including the absolute line $\omega$. This is achieved when
$I_2=0$. Indeed, $a_{11}+a_{22}=0$ and $I_2=a_{11}\cdot a_{22}=0$
entails $a_{11}=a_{22}=0$.\\
The conic section equation (\ref{conich}) turns into
\begin{equation}\label{f}
F(x_0, x_1, x_2)\equiv
2a_{01}x_1x_0+2a_{02}x_2x_0+a_{00}{x_0}^2=x_0(2a_{01}x_1+2a_{02}x_2+a_{00}x_0)=0,
\end{equation}
out of which we read the invariants (\ref{invarijante}):
$$
I_1=0, \quad I_2=0, \quad  I_3=0,  \quad I_4=-{a_{01}}^2+
{a_{02}}^2, \quad I_5=a_{00}.
$$
According to (\ref{norma}) and (\ref{broj}) the possibilities for
the other line, besides $\omega$ ($x_0=0$), are the following:
\begin{itemize}
\item $I_4>0$ yields the second line in (\ref{f}) is spacelike;
\item $I_4<0$ yields it is a timelike straight line;
\item $I_4=0, \quad a_{01}\neq 0$ reveals the line is isotropic.
\end{itemize}
To end this subsection, for
\begin{itemize}
\item $I_4=0, \quad a_{01}=0, \quad I_5\neq 0$ (\ref{f}) represents a
double absolute line $\omega$;
\item $I_4=0, \quad a_{01}=0, \quad I_5=0$ yields from (\ref{f}) a zero
polynomial.
\end{itemize}
\begin{prop}
In the pseudo-Euclidean plane there are 8 (2 proper + 6 degenerated)
different types of conic sections of the 4th family to distinguish
with respect to the group $B_3$ of motions (see Table 5.).
\end{prop}
We conclude with
\begin{tm}
In the pseudo-Euclidean plane there are 43 (20 proper + 23
degenerated) different types of conic sections to distinguish with
respect to the group $B_3$ of motions (see Table 5.).
\end{tm}

 \vspace{7mm}
\textbf{J. Beban-Brki\'{c}}, Faculty of Geodesy, University of Zagreb,
Ka\v{c}i\'{c}eva $26$, HR-$10 000$ Zagreb, Croatia, e-mail:
jbeban@geof.hr \vspace{3mm}\\
\textbf{M. \v{S}imi\'{c} Horvath}, Faculty of
Architecture, University of Zagreb, Ka\v{c}i\'{c}eva $26$, HR-$10
000$ Zagreb, Croatia, e-mail: marija.simic@arhitekt.hr

\end{document}